\pdfoutput=1
\RequirePackage{ifpdf}
\ifpdf 
\documentclass[pdftex]{sigma}
\else
\documentclass{sigma}
\fi

\newtheorem{TheoremAlph}{Theorem}
\newtheorem{CorollaryAlph}[TheoremAlph]{Corollary}

\numberwithin{equation}{section}

\newtheorem{Theorem}{Theorem}[section]
\newtheorem*{Theorem*}{Theorem}
\newtheorem{Corollary}[Theorem]{Corollary}
\newtheorem{Lemma}[Theorem]{Lemma}

 { \theoremstyle{definition}
\newtheorem{Definition}[Theorem]{Definition}

 }

\usepackage{enumitem}

\begin{document}
\allowdisplaybreaks

\newcommand{\arXivNumber}{2211.14610}

\renewcommand{\PaperNumber}{006}

\FirstPageHeading

\ShortArticleName{Positive Intermediate Ricci Curvature on Fibre Bundles}

\ArticleName{Positive Intermediate Ricci Curvature\\ on Fibre Bundles}

\Author{Philipp REISER~$^{\rm a}$ and David J.~WRAITH~$^{\rm b}$}

\AuthorNameForHeading{P.~Reiser and D.J.~Wraith}

\Address{$^{\rm a)}$~Department of Mathematics, University of Fribourg, Switzerland}
\EmailD{\href{mailto:philipp.reiser@unifr.ch}{philipp.reiser@unifr.ch}}

\Address{$^{\rm b)}$~Department of Mathematics and Statistics, National University of Ireland Maynooth,\\
\hphantom{$^{\rm b)}$}~Maynooth, County Kildare, Ireland}
\EmailD{\href{mailto:david.wraith@mu.ie}{david.wraith@mu.ie}}

\ArticleDates{Received August 15, 2024, in final form January 20, 2025; Published online January 23, 2025}

\Abstract{We prove a canonical variation-type result for submersion metrics with positive intermediate Ricci curvatures. This can then be used in conjunction with surgery techniques to establish the existence of metrics with positive intermediate Ricci curvatures on a wide range of examples which had previously only been known to admit positive Ricci curvature, such as highly connected manifolds and exotic spheres. Further, we extend results of the second author on the moduli space of metrics with positive Ricci curvature to positive intermediate Ricci curvatures.}


\Classification{53C20}

\section{Introduction}

A basic problem in Riemannian geometry is to understand the class of manifolds which satisfy a~given curvature condition, with positive curvature conditions having always played an important role. A fundamental question in this regard is as the following: given a fibre bundle in which both base and fibre (which we will assume to be closed manifolds) admit Riemannian metrics which satisfy a particular curvature condition, when does the total space admit such a~metric? If the curvature condition is positive scalar or positive Ricci curvature, (and assuming the fibre metric is invariant under the action of a Lie structure group), the answer turns out to be ``always''. The construction which establishes the existence of such total space metrics is called the {\it ``canonical variation''}. One begins with a submersion metric on the total space with totally geodesic fibres. This is determined by the given base and fibre metrics, together with a choice of horizontal distribution. One then computes the effect on curvature of shrinking the fibre metrics by a uniform scaling factor. This computation, which rests on the O'Neill formulas for Riemannian submersions, shows that provided the scaling factor is small enough, the total space metric will satisfy the desired curvature condition, see, e.g., \cite[Proposition~9.70]{Be}.

It is clear that one can obtain a wealth of new examples in this way. Indeed, it is hard to overstate the importance of this construction in Riemannian geometry. For example, the first examples of Ricci positive metrics on exotic spheres were established by this method (see~\mbox{\cite{Na,Po}}), and Stolz's celebrated classification result for simply-connected manifolds of positive scalar curvature depends on it in a crucial way \cite{St2}. Moreover, the canonical variation can often be combined with other techniques to produce powerful construction tools (see, for example, \cite{Wr2} and \cite{Re}, as well as examples later in this paper).

Given the centrality of the canonical variation in positive scalar and positive Ricci curvature geometry, it makes sense to ask for which other positive curvature conditions such a result will hold. Note that in the context of positive sectional curvature the analogous statement is not true: the situation here is much more complicated.

The focus of this paper is on positive {\it intermediate Ricci curvatures}. These constitute a~natural family of positive curvature conditions which lie in-between positive Ricci and positive sectional curvature. We begin by recalling the definition of these curvatures.
\begin{Definition}
	Given a point $p$ in a Riemannian manifold $M$, and a collection $v_0,\dots,v_k$ of orthonormal vectors in $T_pM$, the $k^{\rm th}$-intermediate Ricci curvature at $p$ corresponding to this choice of vectors is defined to be $\sum_{i=1}^k K(v_0,v_i)$, where $K$ denotes the sectional curvature.
\end{Definition}
Notice that for an $n$-dimensional manifold, these curvatures interpolate between the Ricci curvature and the sectional curvature as $k$ ranges from $n-1$ down to 1. We will denote the $k^{\rm th}$-intermediate Ricci curvature by ${\rm Ric}_k$.

Our first main result is the following.

\begin{TheoremAlph}\label{main}
	Let $(M,g)$ be a compact Riemannian manifold, and suppose that $(M,g)\xrightarrow{\pi} (B^q,\check{g})$ is a Riemannian submersion with totally geodesic fibres. Suppose that $(B,\check{g})$ has ${\rm Ric}_{k_1}>0$, and that the fibres $\bigl(F_b^p,\hat{g}_b\bigr)=\bigl(\pi^{-1}(b),g|_{\pi^{-1}(b)}\bigr)$ for $b \in B$ $($which are necessarily isometric$)$ have ${\rm Ric}_{k_2}>0$. Then there exists $\tau=\tau(M,g,\pi)>0$ such that for all positive $t<\tau$, scaling $g$ in fibre directions by $t$ yields a metric on $M$ with ${\rm Ric}_k>0$ for every $k\ge \max\{k_1+p,k_2+q\}$.	
\end{TheoremAlph}

Although intermediate curvatures have been studied to some extent for several decades, recent years have seen a significant increase in interest. This perhaps reflects a realisation of their potential to provide a more full view of Riemannian curvature, and to cast the classical measures of curvature in a broader context which could ultimately shed more light on classical problems, such as the geometry of submanifolds (see, e.g., \cite{AQZ,GW1,GW2}), symmetries (see, e.g.,~\cite{Mo1,Mo2}), synthetic notions of classical curvature conditions (see, e.g., \cite{KM1}) and topological implications (see, e.g., \cite{GX,GW3,Sh,Wi,Wu}). Further, examples of manifolds with ${\rm Ric}_k>0$ were constructed in the presence of symmetries in \cite{AQZ,DGM}, and by using surgery techniques in \cite{RW}. Nevertheless, it remains a difficult problem to construct metrics of ${\rm Ric}_k>0$ and examples are still rare, as one of the main difficulties lies in detecting when a metric satisfies such a condition. We will say a~little more about our approach to this problem at the end of this section.

Riemannian submersions with totally geodesic fibres, as demanded by the hypotheses of Theorem \ref{main}, are easy to construct. By a result of Vilms \cite{Vi} (see also \cite[Theorem~9.59]{Be}),
fibre bundles can be equipped with such structures under mild conditions. This immediately yields the following.

\begin{CorollaryAlph}\label{fibre_bundles}
	Let $M\xrightarrow{\pi}B^q$ be a fibre bundle with fibre $F^p$ and structure group $G$ so that base and fibre are compact. Suppose that $B$ admits a metric $\check{g}$ with ${\rm Ric}_{k_1}>0$ and $F$ admits a $G$-invariant metric $\hat{g}$ with ${\rm Ric}_{k_2}>0$. Then $M$ admits a metric with ${\rm Ric}_k>0$ for every $k\geq \max\{k_1+p,k_2+q \}$.
\end{CorollaryAlph}

Corollary \ref{fibre_bundles} is a direct generalization of a result of Gromoll and Walschap for positive Ricci curvature \cite[Theorem 2.7.3]{GW}, see also \cite{CS,Na,Po}.

We observe that Corollary \ref{fibre_bundles} agrees precisely with the situation for Riemannian products, (see \cite[Proposition~2.2]{DGM}). Thus {\it twisting the fibres has no impact on the existence of positive intermediate Ricci curvatures} in this range.

In \cite{RW}, the authors considered Gromov's bound on Betti numbers for manifolds with sectional curvature bounded below~\cite{G}. This was shown to fail for positive Ricci curvature by Sha and Yang \cite{SY}, and in \cite{RW} it was established that this failure continues for positive $k^{\rm th}$ intermediate Ricci curvatures for $k$ down to roughly the middle dimension. This result was proved by establishing a surgery theorem for positive intermediate Ricci curvatures. Taking surgery insights from \cite{RW} and combining them with Corollary~\ref{fibre_bundles} leads to the following.

\begin{TheoremAlph}\label{plumbing}
	Let $W$ be the manifold obtained by plumbing linear disc bundles over spheres according to a simply-connected graph. For a fixed bundle in this graph denote the base dimension by $q+1$ and the fibre dimension by $p+1$ and suppose that $p,q\geq2$. Then $\partial W$ admits a metric with ${\rm Ric}_k>0$ for all $k\geq \max\{p,q \}+2$.
\end{TheoremAlph}

In short, Corollary \ref{fibre_bundles} allows the basic ${\rm Ric}_k>0$ surgery result in \cite{RW} to be enhanced by allowing surgery to be performed using different trivialisations of the normal bundle, which results in Theorem \ref{plumbing}. This is analogous to the improvement of the basic Ricci positive surgery theorem in \cite{SY} achieved by the second author in \cite{Wr2}, and further enhanced by the first author in~\cite{Re}. (For a detailed discussion of surgery with different trivialisations in the context of positive Ricci curvature, see, for example, \cite[Section~3]{Wr3}.)

In \cite{RW}, the authors showed the existence of ${\rm Ric}_k>0$ metrics on the connected sums $\sharp_r S^n \times S^m$ for $n,m \ge 2$ with $n \neq m$, whenever $k \ge \max\{n,m\}+1$, and also on $\sharp_r S^n \times S^n$ for $n \ge 3$ when $k \ge n+2$. Using Theorem~\ref{plumbing}, we can now greatly extend this collection of examples with the following, which generalize the positive Ricci curvature existence results of the second author~\cite{Wr1} in the case of $(1)$ and $(3)$ below, and those of D.~Crowley and the second author~\cite{CW} in the case of~$(2)$.

\begin{TheoremAlph}\label{examples}
	We have the following examples of manifolds with a metric of ${\rm Ric}_k>0$:
	\begin{enumerate}[label=\textnormal{(\arabic*)}]\itemsep=0pt
		\item Every homotopy sphere $\Sigma^{2m-1}$ that bounds a parallelisable manifold admits a metric with ${\rm Ric}_k>0$ for $k\geq m+1$.
		\item Let $M^{4m-1}$, $m\geq2$, be a $(2m-2)$-connected manifold that is $(2m-1)$-parallelisable if $m\equiv1\mod4$. Then there is a homotopy sphere $\Sigma$ so that $M\#\Sigma$ admits a metric with ${\rm Ric}_k>0$ for $k\geq 2m+1$.
		\item The unique exotic $8$-sphere admits a metric with ${\rm Ric}_6>0$, and the unique exotic $16$-sphere admits a metric with ${\rm Ric}_{10}>0$.
	\end{enumerate}
\end{TheoremAlph}

For $m=2,3$ item $(2)$ can be refined: In fact, every 2-connected 7-manifold admits a~metric with ${\rm Ric}_5>0$ and every 4-connected 11-manifold admits a metric with ${\rm Ric}_7>0$, cf.\ \mbox{\cite[Corollary~5.4]{CW}}.

The examples in item $(3)$ are exotic spheres which can be obtained by plumbing and do not bound a parallelisable manifold. Note that while we only listed the most well-known examples of this type, there exist further examples which provide additional applications of Theorem \ref{plumbing}.

For more than a decade, one of the major trends in Riemannian geometry has been the study of spaces and moduli spaces of Riemannian metrics satisfying various forms of positive curvature condition, see, e.g., \cite{TW} and the references therein. In this direction, we consider moduli spaces of positive intermediate Ricci curvature metrics, and propose the following result.

\begin{TheoremAlph}\label{spaces}
	Let $M$ be any of the manifolds appearing in $(1)$ of Theorem \ref{examples} with $m$ even, or appearing in $(2)$ in the following situations: $m$ is odd, or $m$ is even and the Pontryagin class $p_{\frac{m}{2}}(M)\in H^{2m}(M;\mathbb{Z})$ is a torsion element. Then for any $k$ for which Theorem \ref{examples} guarantees us a metric with ${\rm Ric}_k>0$, the moduli space of ${\rm Ric}_k>0$ metrics on $M$ $($in $(1))$, resp.\ $M\#\Sigma$ $($in~$(2))$, has infinitely many path-components.
\end{TheoremAlph}

In \cite{KS}, Kreck and Stolz showed that the positive sectional curvature (i.e., ${\rm Ric}_1>0$) moduli space of certain $7$-dimensional manifolds has at least two path-components. Moreover in \cite{DKT}, for each dimension $4n-1\geq 7$, infinite families of closed manifolds are constructed where the positive Ricci curvature (i.e., ${\rm Ric}_{4n-2}>0$) moduli space has infinitely many path-components, and in \cite{DGM} it is shown that these results in fact hold for ${\rm Ric}_{4n-3}>0$. Besides these results, to the best of our knowledge, Theorem \ref{spaces} is the first result on moduli spaces for curvature conditions stricter than positive Ricci curvature.

To prove Theorem \ref{spaces}, one needs to use the surgery techniques behind Theorem \ref{plumbing} together with ideas from \cite{Wr3}. However establishing this result turns out to require very delicate arguments, and its proof occupies a significant part of the sequel.

It was established independently in \cite{De,Go} that the moduli space of metrics of positive Ricci curvature on every total space $E^7$ of a linear $S^3$-bundle over $S^4$ such that $H^4(E)$ is torsion has infinitely many path-components, and this result has been extended in \cite{We} to all total spaces~$E^{15}$ of linear $S^7$-bundles over $S^8$ such that $H^8(E)$ is torsion. Since the metrics considered in~\cite{We} are submersion metrics as in Corollary \ref{fibre_bundles} (or, more precisely, as in Corollary \ref{fibre_bundles1} below), we directly obtain that these results in fact hold for ${\rm Ric}_5>0$ and ${\rm Ric}_9>0$, respectively. We note that this fact also follows from Theorem \ref{spaces} since all these spaces are contained in $(2)$ of Theorem \ref{examples}, and it follows from the proof of Theorem \ref{examples} that for these spaces the corresponding homotopy sphere appearing in Theorem \ref{examples} can be chosen to be a standard sphere.

In order to establish Theorem \ref{main} and Corollary \ref{fibre_bundles}, we study the following operator. For $X \in T_pM$, define
\begin{align*}
	R_X\colon \ T_pM &\to T_pM, \\
	Y &\mapsto R(X,Y)X
\end{align*}
(Note that, depending on the chosen sign convention for the curvature tensor, there exist different conventions to define this operator. See Section \ref{formulas} for our choice of sign convention.)
It is easy to see that this operator is symmetric, so there exists an orthonormal basis of eigenvectors. Given orthonormal unit vectors $X$, $Y$, if $Y$ is an eigenvector of $R_X$ with eigenvalue $\lambda$, then \[\lambda=\langle R_X(Y),Y\rangle=K(X,Y).\] Thus we can obtain $k^{\rm th}$ intermediate Ricci curvatures by summing eigenvalues of $R_X$. More precisely, we have the following.
\begin{Lemma}\label{k+1}
	A Riemannian manifold $M$ has ${\rm Ric}_k>0$ if and only if for every unit vector $X \in TM$, the sum of any $(k+1)$-eigenvalues of $R_X$ is positive.
\end{Lemma}
Note that this is essentially \cite[Lemma 2.6]{RW}, and we refer the reader there for more details.~(Of course $X$ is an eigenvector of $R_X$ with eigenvalue~0, and therefore $\langle R_X(X),X\rangle =K(X,X)=0$, hence the need to sum $k+1$ eigenvectors to obtain a positive sum.) Theorem \ref{main} is then obtained by a careful analysis of the eigenvectors of $R_X$.

This paper is laid out as follows: in Section \ref{formulas}, we examine the canonical variation in the context of intermediate Ricci curvatures, and in Section \ref{proofs}, we establish all the main theorems.

\section{Curvature formulae}\label{formulas}

In this section, we present some formulas for the curvature tensor of a Riemannian submersion. We are interested in the case where the fibres are totally geodesic, and the effect the canonical variation has on the curvature tensor in this case. Recall that the canonical variation involves scaling the fibre metrics by a uniform constant, leaving the base metric and the horizontal distribution untouched. The case of particular interest is when the fibre scaling factor approaches zero. Note that the canonical variation preserves the property of being totally geodesic.

The starting point for the curvature formulas we will require is the O'Neill formulas for Riemannian submersions, as presented, for example, in~\cite{Be}. As this book is the definitive reference for the curvature of Riemannian submersions and the canonical variation, we will adopt the conventions used there. In particular, note that for the canonical variation, we will scale in fibre directions by a constant $t>0$ \big(as opposed to the more usual $t^2$\big). We also use the following curvature tensor convention:
\[
\langle R(V,W)X,Y\rangle=\bigl\langle \nabla_W\nabla_V X -\nabla_V \nabla_W X+\nabla_{[V,W]}X,Y\bigr\rangle,
\]
where $\langle\,,\,\rangle$ denotes the Riemannian metric.

The formulas below are due to O'Neill and Gray. They appear in \cite{Be} as Theorem 9.28. In~stating these, we have specialized to the case of totally geodesic fibres, i.e., we have assumed that the $T$-tensor is identically zero. Recall that the \emph{$A$-tensor} of a Riemannian submersion is a~$(2,1)$-tensor field which vanishes if and only if the horizontal distribution is integrable. We refer to \cite[Section 9.C]{Be} for its definition and basic properties.

\begin{Lemma}\label{formulas1}
	Given a Riemannian submersion with totally geodesic fibres, let $X$, $Y$, $Z$, $Z'$ denote arbitrary horizontal vectors, and $U$, $V$, $W$, $W'$ arbitrary vertical vectors, all over the same base point. Then
	\begin{align*}
		&\bigl\langle R(U,V)W,W' \bigr\rangle =\bigl\langle \hat{R}(U,V)W,W' \bigr\rangle, \\
		&\langle R(U,V)W,X \rangle =0, \\
		&\langle R(X,U)Y,V \rangle =\langle (\nabla_U A)_X Y,V \rangle+\langle A_X U, A_Y V\rangle, \\
		&\langle R(U,V)X,Y \rangle =\langle (\nabla_U A)_X Y,V \rangle-\langle (\nabla_V A)_X Y,U\rangle+\langle A_XU,A_YV \rangle-\langle A_XV, A_YU \rangle, \\
		&\langle R(X,Y)Z,U \rangle = \langle (\nabla_ZA)_X Y,U \rangle, \\
		&\bigl\langle R(X,Y)Z,Z' \bigr\rangle =\bigl\langle \check{R}(X,Y)Z,Z' \bigr\rangle-2\bigl\langle A_XY,A_ZZ' \bigr\rangle + \bigl\langle A_Y Z,A_X Z'\bigr\rangle-\bigl\langle A_X Z, A_Y Z'\bigr\rangle.
	\end{align*}
	Here, $\hat{R}$ is the curvature tensor of the fibre metric, and $\check{R}$ that of the base metric.
\end{Lemma}

The next lemma is an adaptation of \cite[Lemma 9.69]{Be}, where we only state the effect of the canonical variation on formulas which appear in Lemma \ref{formulas1}.

\begin{Lemma}\label{formulas2}
	Given a Riemannian submersion with fibre metric $\hat{g}$, we scale in fibre directions by $t>0$, that is, we take the new fibre metric to be $t\hat{g}$, but leave the base metric and horizontal distribution unchanged. Denoting quantities for the rescaled metric with a superscript $t$, and using the same vector notation as in Lemma~$\ref{formulas1}$, the following formulas hold:
	\begin{align*}
		&A^t_X Y=A_XY, \\
		&A^t_X U=tA_XU, \\
		&[\langle (\nabla_X A)_YZ,U\rangle]^t=t\langle (\nabla_XA)_Y Z,U\rangle, \\
		&[\langle (\nabla_UA)_X Y,V\rangle]^t=t\langle (\nabla_U A)_X Y,V \rangle+\bigl(t-t^2\bigr)\bigl(\langle A_XU,A_YV\rangle-\langle A_XV,A_YU\rangle\bigr).
	\end{align*}
\end{Lemma}

It can be deduced from the Koszul formula, for example, that scaling a metric $g$ by a constant factor, i.e., replacing $g$ by $t\cdot g$ for $t>0$, leaves the Levi-Civita connection unchanged. Therefore given vectors $U$, $V$, $W$, we see immediately that $R(U,V)W$ is unaffected by scaling the metric as well. Hence, in the case of arbitrary vertical vectors $U$, $V$, $W$, $W'$, we then have
\[
\bigl[\bigl\langle\hat{R}(U,V)W,W'\bigr\rangle\bigr]^t=t\bigl\langle\hat{R}(U,V)W,W'\bigr\rangle.
\]
Further, since the metric is unchanged in horizontal directions, for horizontal vectors~$X$,~$Y$,~$Z$,~$Z'$, we have
\[
\bigl[\bigl\langle \check{R}(X,Y)Z,Z'\bigr\rangle\bigr]^t=\bigl\langle \check{R}(X,Y)Z,Z'\bigr\rangle.
\]
Combining these observations with Lemma \ref{formulas2}, a computation yields the following formulas, which show how the expressions Lemma \ref{formulas1} transform under the canonical variation.

\begin{Corollary}\label{formulas3}
	Under the same assumptions and with the same notation as the previous two lemmas, the following formulas hold:
	\begin{align*}
		&\bigl[\bigl\langle R(U,V)W,W' \bigr\rangle\bigr]^t =t\bigl\langle \hat{R}(U,V)W,W' \bigr\rangle, \\
		&[\langle R(U,V)W,X\rangle]^t =0, \\
		&[\langle R(X,U)Y,V \rangle]^t =t\langle (\nabla_U A)_X Y,V \rangle+\bigl(t-t^2\bigr)(\langle A_XU,A_YV\rangle-\langle A_XV,A_YU\rangle)\\
		&\hphantom{[\langle R(X,U)Y,V \rangle]^t =}{} +t^2\langle A_X U, A_Y V\rangle \\
		&\hphantom{[\langle R(X,U)Y,V \rangle]^t }{} =t\langle (\nabla_U A)_X Y,V \rangle +t\langle A_X U, A_Y V\rangle+\bigl(t^2-t\bigr)\langle A_XV,A_YU\rangle,\\
		&[\langle R(U,V)X,Y \rangle]^t =t\langle (\nabla_U A)_X Y,V \rangle-t\langle (\nabla_V A)_X Y,U\rangle+t^2\langle A_XU,A_YV \rangle-t^2\langle A_XV, A_YU \rangle \\
		&\hphantom{[\langle R(U,V)X,Y \rangle]^t =}{} +\bigl(t-t^2\bigr)(\langle A_XU,A_YV\rangle-\langle A_XV,A_YU\rangle\\
&\hphantom{[\langle R(U,V)X,Y \rangle]^t =}{}- \langle A_XV,A_YU\rangle + \langle A_XU,A_YV\rangle)\\
		&\hphantom{[\langle R(U,V)X,Y \rangle]^t }{} =t\langle (\nabla_U A)_X Y,V \rangle-t\langle (\nabla_V A)_X Y,U\rangle\\
		&\hphantom{[\langle R(U,V)X,Y \rangle]^t =}{} +\bigl(2t-t^2\bigr)(\langle A_XU,A_YV\rangle-\langle A_XV,A_YU\rangle),\\
		&[\langle R(X,Y)Z,U \rangle]^t = t\langle (\nabla_ZA)_X Y,U \rangle, \\
		&\bigl[\bigl\langle R(X,Y)Z,Z' \bigr\rangle\bigr]^t=\bigl\langle \check{R}(X,Y)Z,Z' \bigr\rangle\!-2t\bigl\langle A_XY,A_ZZ' \bigr\rangle\! + t\bigl\langle A_Y Z,A_X Z'\bigr\rangle\!-t\bigl\langle A_X Z, A_Y Z'\bigr\rangle.
	\end{align*}
\end{Corollary}

Let us fix the following notation. Let $g$ denote the original (unscaled) Riemannian submersion metric, and denote by $g_t$ the canonical variation metric corresponding to the vertical scaling by~$t$. We now wish to specialize to the case where all vectors are unit vectors for $g_t$. Suppose that~$U$,~$V$,~$W$,~$W'$ are all unit vertical vectors for $g$ (i.e., the case $t=1$), then the corresponding unit vectors for $g_t$ are $U/\sqrt t$, $V/\sqrt t$, $W/\sqrt t$, $W'/\sqrt t$. From Corollary \ref{formulas3}, we immediately obtain the expressions below.

\begin{Corollary}\label{formulas4}
	Suppose that $X$, $Y$, $Z$, $Z'$ are unit horizontal vectors for $g_t$. If $U$, $V$, $W$, $W'$ are all unit vertical vectors for $g$, and $U_t$, $V_t$, $W_t$, $W'_t$ the corresponding unit vectors for $g_t$, we~have
	\begin{align*}
		&\bigl[\bigl\langle R(U_t,V_t)W_t,W'_t \bigr\rangle\bigr]^t =\frac{1}{t}\bigl\langle \hat{R}(U,V)W,W' \bigr\rangle, \\
		&[\langle R(U_t,V_t)W_t,X\rangle]^t =0, \\
		&[\langle R(X,U_t)Y,V_t \rangle]^t =\langle (\nabla_U A)_X Y,V \rangle +\langle A_X U, A_Y V\rangle+(t-1)\langle A_XV,A_YU\rangle,\\
		&[\langle R(U_t,V_t)X,Y \rangle]^t =\langle (\nabla_U A)_X Y,V \rangle-\langle (\nabla_V A)_X Y,U\rangle\\
		&\hphantom{[\langle R(U_t,V_t)X,Y \rangle]^t =}{} +(2-t)(\langle A_XU,A_YV\rangle-\langle A_XV,A_YU\rangle),\\
		&[\langle R(X,Y)Z,U_t \rangle]^t = \sqrt{t}\langle (\nabla_ZA)_X Y,U \rangle, \\
		&\bigl[\bigl\langle R(X,Y)Z,Z' \bigr\rangle\bigr]^t =\bigl\langle \check{R}(X,Y)Z,Z' \bigr\rangle\!-2t\bigl\langle A_XY,A_ZZ' \bigr\rangle\! + t\bigl\langle A_Y Z,A_X Z'\bigr\rangle\!-t\bigl\langle A_X Z, A_Y Z'\bigr\rangle.
	\end{align*}
\end{Corollary}

\begin{Corollary}\label{R_matrix}
	Let $\Psi$ be a unit tangent vector with respect to $g_t$, i.e., we can write $\Psi=\lambda U_t+\mu X$ with $\lambda,\mu\in\mathbb{R}$ so that $\lambda^2+\mu^2=1$. The map $R_\Psi^t$, according to the splitting into vertical and horizontal parts, is given by the matrix
	\begin{equation*}
	\begin{pmatrix}
		\frac{\lambda^2}{t}\hat{R}_U+O(t) & O(|\lambda|)+O\bigl(\sqrt{t}\bigr)\\
		O(|\lambda|)+O\bigl(\sqrt{t}\bigr) & \mu^2 \check{R}_X+O(|\lambda|)+O(t)
	\end{pmatrix}
\end{equation*}
	as $t,\lambda\to0$.
\end{Corollary}
Here for a linear map $T_\varepsilon$ that depends on a parameter $\varepsilon$ we say that $T_\varepsilon=O(f(\varepsilon))$ for a real function $f$, if the operator norm $\lVert T_\varepsilon\rVert$ satisfies $\lVert T_\varepsilon\rVert=O(f(\varepsilon))$ in the usual sense.
\begin{proof}
	We first consider the first entry of the matrix. Here, we have
	\begin{align*}
		[\langle R(\Psi,V_t)\Psi,W_t\rangle]^t&{}=\lambda^2[\langle R(U_t,V_t)U_t,W_t\rangle]^t+\mu^2 [\langle R(X,V_t)X,W_t\rangle]^t\\
		&\quad{} + \lambda\mu\bigl([\langle R(U_t,V_t)X,W_t\rangle]^t + [\langle R(X,V_t)U_t,W_t\rangle]^t \bigr).
	\end{align*}
	By Corollary \ref{formulas4}, we have
	\begin{align*}
		[\langle R(U_t,V_t)U_t,W_t\rangle]^t &{}= \frac{1}{t}\bigl\langle \hat{R}(U,V)U,W\bigr\rangle,\\
		[\langle R(X,V_t)X,W_t\rangle]^t &{}= \langle (\nabla_V A)_X X,W\rangle +t\langle A_X V,A_X W\rangle,\\
		[\langle R(U_t,V_t)X,W_t\rangle]^t &{}=0,\\
		[\langle R(X,V_t)U_t,W_t\rangle]^t &{}=0.
	\end{align*}
	It can be checked that $\nabla_V A$ enjoys the same antisymmetry property with respect to pairs of horizontal vectors as the $A$-tensor itself. Thus $\langle (\nabla_V A)_X X,W\rangle=0$. Since the $A$-tensor is bounded (as is any operator on a finite-dimensional vector space) and since $|\mu|\leq 1$, the first entry of the matrix follows.
	
	Similarly, for the second entry we have
	\begin{align*}
		[\langle R(\Psi,V_t)\Psi,Y\rangle ]^t&{}= \lambda^2[\langle R(U_t,V_t)U_t,Y\rangle]^t+\mu^2 [\langle R(X,V_t)X,Y\rangle]^t\\
		&\quad{} + \lambda\mu\bigl([\langle R(U_t,V_t)X,Y\rangle]^t + [\langle R(X,V_t)U_t,Y\rangle]^t \bigr).
	\end{align*}
	Again, by Corollary \ref{formulas4}, we have
	\begin{align*}
		[\langle R(U_t,V_t)U_t,Y\rangle]^t&{}=0,\\
		[\langle R(X,V_t)X,Y\rangle]^t&{}=\sqrt{t}\langle(\nabla_X A)_X Y,V\rangle=O\bigl(\sqrt{t}\bigr),
	\end{align*}
	and both $[\langle R(U_t,V_t)X,Y\rangle]^t$ and $[\langle R(X,V_t)U_t,Y\rangle]^t$ are bounded as $t\to 0$. Hence, we obtain, again by using $|\mu|\leq 1$,
	\begin{align*}
		[\langle R(U_t,V_t)U_t,Y\rangle]^t&{}=O\bigl(\sqrt{t}\bigr)+O(|\lambda|).
	\end{align*}
	By the symmetry properties of $R_\Psi$, we obtain the same conclusion for the third entry.
	
	Finally, for the fourth entry, we have
	\begin{align*}
		[\langle R(\Psi,Y)\Psi,Z\rangle ]^t&{}= \lambda^2[\langle R(U_t,Y)U_t,Z\rangle]^t+\mu^2 [\langle R(X,Y)X,Z\rangle]^t\\
		&\quad{}+ \lambda\mu\bigl([\langle R(U_t,Y)X,Z\rangle]^t + [\langle R(X,Y)U_t,Z\rangle]^t \bigr)
	\end{align*}		
	By Corollary \ref{formulas4}, the terms $[\langle R(U_t,Y)U_t,Z\rangle]^t$, $[\langle R(U_t,Y)X,Z\rangle]^t$, and $[\langle R(X,Y)U_t,Z\rangle]^t$ are bounded as $t\to 0$ and
	\[ [\langle R(X,Y)X,Z\rangle]^t=\bigl\langle\check{R}(X,Y)X,Z\bigr\rangle+O(t). \]
	Hence, using $|\mu|,|\lambda|\leq 1$, we obtain
	\begin{align*}
		[\langle R(\Psi,Y)\Psi,Z\rangle ]^t&{}= O(|\lambda|)+\mu^2\bigl(\bigl\langle \check{R}(X,Y)X,Z\bigr\rangle+O(t) \bigr)\\
		&{}= \mu^2\bigl(\bigl\langle \check{R}(X,Y)X,Z\bigr\rangle\bigr)+O(|\lambda|)+O(t).
 \tag*{\qed}
 \end{align*}
 \renewcommand{\qed}{}
\end{proof}

\section{Proof of the main theorems}\label{proofs}

\subsection{Proof of Theorems \ref{main}--\ref{examples}}

\begin{proof}[Proof of Theorem \ref{main}]
	We fix a point $p\in M$. As in Corollary \ref{R_matrix}, any unit tangent vector $\Psi\in T_p M$ with respect to $g_t$ can be written as $\Psi=\lambda U_t+\mu X$ with $\lambda^2+\mu^2=1$ and $\lambda,\mu\geq 0$. In~the~following we therefore fix $U$ and $X$, consider $\Psi=\lambda U_t+\mu X$ as parametrized by $\lambda $ and $t$, and show that for all $t$ sufficiently small and all $\lambda\in[0,1]$, the sum of any $(k+1)$ eigenvalues of the map $R^t_\Psi$ is positive for $k\geq\max\{k_1+p,k_2+q \}$. (Recall that by Lemma~\ref{k+1}, this is equivalent to the condition ${\rm Ric}_k>0$.)
	
	By Corollary \ref{R_matrix}, the map $R_\Psi^t$ is given by the matrix
	\[
	\begin{pmatrix}
		\frac{\lambda^2}{t}\hat{R}_U+O(t) & O(|\lambda|)+O\bigl(\sqrt{t}\bigr)\\
		O(|\lambda|)+O\bigl(\sqrt{t}\bigr) & \mu^2 \check{R}_X+O(|\lambda|)+O(t)
	\end{pmatrix}
	\]
	as $t,\lambda\to0$. For $\lambda=0$, we obtain the matrix
	\[
	\begin{pmatrix}
		O(t) & O\bigl(\sqrt{t}\bigr)\\
		O\bigl(\sqrt{t}\bigr) & \check{R}_X+O(t)
	\end{pmatrix}.
	\]
	By assumption $\check{R}_X$ is $k_1$-positive, so the sum of any $k_1+1$ of its eigenvalues is positive. Since eigenvalues depend continuously on the matrix, see, e.g., \cite[Corollary VI.1.6]{Bh}, there exists $t_0>0$ so that the sum of any $p+k_1+1$ eigenvalues of the above matrix is positive for all $t\in(0,t_0)$. Since $R^t_\Psi$ depends continuously on $\Psi$, this conclusion holds for all $\lambda\in[0,\lambda_0)$ for some $\lambda_0>0$.
	
	Now let $\lambda\geq\lambda_0$. Multiplying the matrix by $t$ gives the matrix
	\[
	\begin{pmatrix}
		\lambda^2 \hat{R}_U+O\bigl(t^2\bigr) & t\cdot O(|\lambda|)+O\bigl(t^{\frac{3}{2} }\bigr)\\
		t\cdot O(|\lambda|)+O\bigl(t^{\frac{3}{2}}\bigr) & \mu^2t \check{R}_X+t\cdot O(|\lambda|)+O\bigl(t^2\bigr)
	\end{pmatrix},
	\]
	and note that the condition of $k$-positivity is invariant under scalar multiplication of the matrix by a positive constant. As the fibre is compact, $\hat{R}_U$ is $(k_2+1)$-bounded away from zero. Thus since $\lambda\in[\lambda_0,1]$, and the $O(|\lambda|)$-terms depend continuously on $\lambda$ and hence are bounded for $\lambda\in[\lambda_0,1]$, there exists $t_1>0$ so that the sum of any $k_2+q+1$ eigenvalues of this matrix is positive for all $t\in(0,t_1)$. 
	
	We set $\tau=\min(t_0,t_1)$, so the sum of any $k+1$ eigenvalues of $R_\Psi^t$ is positive for all $t\in(0,\tau)$ and $\lambda\in[0,1]$, and for all $k\geq \max\{k_1+p,k_2+q\}$. By continuity and compactness of the space of unit horizontal respectively vertical vectors, it then follows that there exists a uniform choice of $t$ for all unit tangent vectors $\Psi$, and by compactness of $M$, for all $p\in M$.
\end{proof}

Corollary \ref{fibre_bundles} is now a special case of the following corollary.

\begin{Corollary}\label{fibre_bundles1}
	Let $M\xrightarrow{\pi}B^q$ be a fibre bundle with fibre $F^p$ and structure group $G$ so that base and fibre are compact. Suppose that $B$ admits a metric $\check{g}$ with ${\rm Ric}_{k_1}>0$ and $F$ admits a $G$-invariant metric $\hat{g}$ with ${\rm Ric}_{k_2}>0$. Then for any horizontal distribution $\mathcal{H}$ on $\pi$ that is obtained from a principal connection on the associated principal $G$-bundle, and for all sufficiently small $t>0$, $M$ admits a metric $g$ with ${\rm Ric}_k>0$ for every $k\ge \max\{k_1+p,k_2+q\}$, so that $(M,g)\to (B,\check{g})$ is a Riemannian submersion with totally geodesic fibres isometric to $(F,t\hat{g})$ and horizontal distribution $\mathcal{H}$.
\end{Corollary}

\begin{proof}
	All claims in Corollary \ref{fibre_bundles1} are a direct consequence of the so-called Vilms construction, see, e.g., \cite[Theorem~9.59]{Be}, except the claim about ${\rm Ric}_k>0$. This claim directly follows from Theorem~\ref{main}.
\end{proof}

For the proof of Theorem \ref{plumbing}, we need the theorem below, which extends \cite[Theorem~C]{RW}. For~${\rho>0}$ we denote by $S^p(\rho)$ the round sphere of radius $\rho$ and for $R,N>0$ we denote by~\smash{$D^{q+1}_R(N)$} a~geodesic ball of radius $R$ in $S^{q+1}(N)$. By analogy with \cite{Wr2} we make the following assumptions:
\begin{enumerate}[label=\textnormal{(\arabic*)}]\itemsep=0pt
	\item $(M^n,g_M)$ is a Riemannian manifold with ${\rm Ric}_k>0$.
	\item $\iota\colon S^p(\rho)\times D_R^{q+1}(N)\hookrightarrow (M,g_M)$, $p+q+1=n$, is an isometric embedding.
	\item $T\colon S^p\to {\rm SO}(q+1)$ is a smooth map. Define the diffeomorphism $\tilde{T}\colon S^p\times D^{q+1}\to S^p\times D^{q+1},\, (x,y)\mapsto (x,T_x(y))$.
\end{enumerate}

\begin{Theorem}\label{surgery}
	Under the assumptions $(1)$--$(3)$, suppose $k\geq\max\{p,q \}+2$ and $p,q\geq2$. Then there exists a constant $\kappa=\kappa(p,q,k,R/N,T)>0$ such that if $\frac{\rho}{N}<\kappa$, then the manifold
	\[M_{\tilde{T}\circ\iota}=M\setminus\mathrm{im}(\iota)^\circ\cup_{\tilde{T}|_\partial} \bigl(D^{p+1}\times S^q\bigr) \]
	admits a metric with ${\rm Ric}_k>0$. Further, for any $\rho',R',N'>0$ with $\frac{R'}{N'}<\frac{\pi}{2}$ and $\rho'$ sufficiently small we can construct the metric to coincide near the centre of $D^{p+1}\times S^q$ with $D^{p+1}_{R'}(N')\times S^q(\rho')$. In this case, the constant $\kappa$ additionally depends on $\rho'$, $R'$, $N'$.
\end{Theorem}
The proof is an adaptation of the proof of \cite[Theorem 0.3]{Wr2} to our situation.
\begin{proof}
	Let $R''\in\bigl(R',N'\frac{\pi}{2}\bigr)$. By \cite[Theorem C and Remark 4.2]{RW}, there is a metric on $D^{p+1}\times S^q$ with ${\rm Ric}_k>0$ for $k\geq\max\{p,q\}+2$ that glues smoothly with $M\setminus\mathrm{im}(\iota)^\circ$, and near the centre of~${D^{p+1}\times S^q}$, after possibly rescaling, takes the form $D^{p+1}_{R''}(N')\times S^q(\rho')$. We will now replace this part by a linear sphere bundle equipped with a submersion metric so that the twisting described by $T$ is taken into account.
	
	Let $E\xrightarrow{\pi}S^{p+1}$ be the linear $S^q$-bundle, whose clutching function is given by $T$. We equip the base $S^{p+1}$ with the metric $S^{p+1}(N')$, which can be written as
	\[{\rm d}t^2+{N'}^2\sin^2\biggl(\frac{t}{N'} \biggr){\rm d}s_p^2. \]	
	Here we identified $S^{p+1}$ with the space obtained from $[0,N'\pi]\times S^p$ by collapsing each of the boundary components $\{0 \}\times S^p$ and $\{N'\pi \}\times S^p$.
	
	Now let $\mathcal{H}$ be a horizontal distribution on $\pi$ that is integrable in a neighbourhood of both $[0,R']\times S^p$ and $[R'',N'\pi]\times S^p$. Note that such a distribution exists since $\pi$ restricted to each of these neighbourhoods is trivial, as the bases of these restricted bundles are discs, and hence contractible. For $\rho'$ sufficiently small, Corollary \ref{fibre_bundles1} now provides a submersion metric $g_E$ on $E$ with totally geodesic fibres isometric to $\rho'{\rm d}s_q^2$ and horizontal distribution $\mathcal{H}$ that has ${\rm Ric}_k>0$ for $k\geq\max\{p+1,q+1 \}$. Since $\mathcal{H}$ is integrable on $[0,R']\times S^p$ and $[R'',N'\pi]\times S^p$, the metric $g_E$ is isometric to a product metric over these parts. In particular, it is isometric to $D^{p+1}_{R'}(N')\times S^q(\rho')$ over $[0,R']\times S^p$ and, when restricted to the part of the bundle over $[0,R'']\times S^p$, a neighbourhood of the resulting boundary is isometric to a neighbourhood of the boundary in $D^{p+1}_{R''}(N')\times S^q(\rho')$.
	
Denoting $\pi^{-1}([0,R'']\times S^p)$ by $E'$, we can replace the part of our original metric on $D^{p+1}\times S^q$ which is isometric to $D^{p+1}_{R''}(N')\times S^q(\rho')$ by $E'$ equipped with the restriction of the metric $g_E$, and hence obtain a smooth metric with ${\rm Ric}_k>0$ for $k\geq \max\{p,q\}+2$. The resulting manifold is precisely $M_{\tilde{T}\circ \iota}$.
\end{proof}

\begin{proof}[Proof of Theorem \ref{plumbing}]
	The proof goes along the same lines as the proof of the corresponding statement for positive Ricci curvature \cite[Theorem 2.2]{Wr1} and \cite[Theorem B]{Re} with the surgery theorems used in each case replaced by Theorem \ref{surgery}. For convenience, we recall the argument below.
	
	First observe that the effect on the boundary of plumbing two linear disc bundles over spheres is precisely a surgery operation as in Theorem \ref{surgery}, with $M$ being the sphere bundle of the first bundle, $\iota$ the embedding of a local trivialization of a fibre sphere and $T$ the clutching function of the second disc bundle.
	
	The metric on $M$ is a submersion metric with totally geodesic and round fibres, according to a horizontal distribution that is flat over the embedded disc, so that the metric is a product metric on this part. By Corollary \ref{fibre_bundles1}, this metric has ${\rm Ric}_k>0$ for all $k\geq\max\{p,q \}+2$ for sufficiently small fibre radius.
	
	We now pick a root of the graph and equip the corresponding bundle with a submersion metric with totally geodesic and round fibres according to a horizontal distribution that is integrable over a finite number of embedded discs, where each disc corresponds to a connecting edge. As~a~consequence, the metric is a product metric over these parts. By Corollary \ref{fibre_bundles1}, this metric has ${\rm Ric}_k>0$ for all $k\geq\max\{p,q \}+2$ for sufficiently small fibre radius. We now apply Theorem \ref{surgery} for each edge connected to the root, where in each case the map $T$ is the clutching function of the bundle connected via the corresponding edge. The inequality $\frac{\rho}{N}<\kappa$ can then be satisfied by possibly choosing a smaller fibre radius for the bundle at the root.
	
Since on the attached parts, we again have isometric embeddings of the form $D^{p+1}_{R'}(N')\times S^q(\rho')$, we can repeat this argument and attach all bundles that have distance~2 from the root. The assumptions in Theorem \ref{surgery} can then be satisfied by possibly reducing the fibre radii of all preceding bundles. In this way we can construct a metric with ${\rm Ric}_k>0$ on the whole manifold~$\partial W$. This concludes the proof.
\end{proof}

\begin{proof}[Proof of Theorem \ref{examples}]
By \cite{Ke, KM}, \cite[Proposition 1.5]{Wr1} and \cite[Theorem C]{CW}, all manifolds considered in $(1)$ and $(2)$ of Theorem \ref{examples} can be obtained as boundaries of plumbings of linear disc bundles over spheres according to a simply-connected graph, where fibre and base dimensions are all equal. Further, the unique exotic $8i$-sphere, $i=1,2$, can be obtained as the boundary of plumbing of a linear $D^{4i+1}$-bundle over $S^{4i}$ with a linear $D^{4i}$-bundle over $S^{4i+1}$, see, e.g.,~\mbox{\cite[Theorem 3 and subsequent Example]{Fr}} for $i=1$, or \cite[Satz~12.1]{St}. The claim now follows directly from Theorem~\ref{plumbing}.
\end{proof}

\subsection{Proof of Theorem \ref{spaces}}

For the proof of Theorem \ref{spaces}, the key arguments we need are exactly the same as those developed for positive Ricci curvature moduli spaces in \cite{Wr3}, specifically the metric deformation arguments in Section 4 of that paper, the extension results in Section 5, and the boundary adjustment and $s$-invariant calculations in Section 6.

Note first that the techniques in \cite{Wr3} are only applied to homotopy spheres which bound parallelisable manifolds, i.e., to the examples listed in $(1)$ of Theorem \ref{examples}. However, the techniques are also applicable to more general boundaries of plumbings, and in particular apply to many of the highly connected manifolds listed in $(2)$ of Theorem \ref{examples}. To be precise, \cite[Theorem 6.2.3.3]{TW} asserts that the arguments in \cite{Wr3} show that moduli space of positive Ricci curvature metrics has infinitely many path-components for each of the manifolds in the statement of Theorem \ref{spaces}.

All we need for the technical results in \cite{Wr3} to work in our current ${\rm Ric}_k>0$ setting is that the scaling functions in the doubly warped product metric used in the surgery process behind Theorem \ref{plumbing}, i.e., in Theorem \ref{surgery}, satisfy a number of basic properties. The properties in question are listed in \cite[Theorem 3.2]{Wr3}, and reflect key details of the double warped product metric constructed in \cite{Wr2}. We have listed these properties below as \ref{EQ:SCAL_FCTS1}--\ref{EQ:SCAL_FCTS9}. The specific definition of the scaling functions is therefore not important.

The main difference between our current surgery process compared to that used in \cite{Wr2} to construct the (Ricci positive) examples studied in \cite{Wr3}, is that the scaling functions are not quite the same, though they share the same basic features. The relevant properties the scaling functions need to satisfy are listed in \cite[Theorem 3.2]{Wr3} and are given as follows. For given $\Delta\in(0,1)$ there is a constant $\rho_0=\rho_0(p,q,\Delta)\in(0,1)$ and a constant $R'\in(0,1/4)$, such that for all $\rho'\leq \rho_0$ there exists $a=a(p,q,\Delta,\rho')$ and smooth functions $h,f\colon [0,a]\to \mathbb{R}_{\geq0}$, so that both $a$ and the functions $h$ and $f$ depend smoothly on $\rho'$, with $a\to\infty$ as $\rho'\to0$, and such that the following properties are satisfied:
\begin{enumerate}[label=\textnormal{(S\arabic*)}]\itemsep=0pt
	\item \label{EQ:SCAL_FCTS1} $f\equiv \rho'$ on $[0,R']$,
	\item \label{EQ:SCAL_FCTS2} $f>0$, $f'\geq 0$ and $f''\geq 0$,
	\item \label{EQ:SCAL_FCTS3} $f'\equiv \Delta$ in a small neighbourhood of $r=a$,
	\item \label{EQ:SCAL_FCTS4} $h(r)=\sin(r)$ for all $r\in[0,R']$,
	\item \label{EQ:SCAL_FCTS5} $h$ is independent of all parameters on $[0,1/2]$,
	\item \label{EQ:SCAL_FCTS6} $h>0$ on $(0,a]$, $h'\geq0$ and $h''\leq 0$,
	\item \label{EQ:SCAL_FCTS7} $h'\equiv 0$ in a neighbourhood of $r=a$,
	\item \label{EQ:SCAL_FCTS8} $\sup_{r\in[0,a]}|h''(r)|=\sup_{r\in[0,1/2]}|h''(r)|$,
	\item \label{EQ:SCAL_FCTS9} $h(a)/f(a)\to0$ as $\rho'\to0$.
\end{enumerate}
We also require that the resulting metric on $D^{p+1}\times S^q$ has non-negative Ricci curvature and positive scalar curvature, and strictly positive Ricci curvature for $r$ small. This last condition on the curvature gets replaced in our setting by the condition ${\rm Ric}_k>0$.

Note that a choice of connection also appears in \cite[Theorem 3.2]{Wr3}. It turns out that the influence of this on the constructions we are considering can always be rendered insignificant by choosing $\rho_0$ smaller if necessary. Thus we can omit the connection from our considerations.

It is easily checked that the functions in \cite{RW}, which are used in the proof of Theorem \ref{surgery} and are based on the functions defined in \cite{Re}, satisfy all of the properties \ref{EQ:SCAL_FCTS1}--\ref{EQ:SCAL_FCTS7}, except for property \ref{EQ:SCAL_FCTS3}. We will see below that property \ref{EQ:SCAL_FCTS3} is not essential for the arguments in \cite{Wr3} to work, and further that property \ref{EQ:SCAL_FCTS9} also holds. However, property \ref{EQ:SCAL_FCTS8} does not hold since for the construction in \cite{RW}, in order to satisfy the required boundary conditions, the function $h$ is bent with very large second derivative. Thus, we need to modify the construction of the scaling functions, which we will now describe.

We follow the first part of \cite[proof of Theorem 4.1]{RW} and apply \cite[Remark 4.2]{RW}, that is, we obtain $t_0<t_2$ and smooth functions $h$, $f$ defined on $[t_0,t_2]$ such that
\begin{align*}
	&h(t)= N'' \sin\biggl(\frac{t-t_0}{ N''}\biggr)\qquad \text{for all}\quad t\in\bigl[t_0,t_0+\alpha R''\bigr],\\
	&f(t)=\rho''\qquad \text{for all}\quad t\in\bigl[t_0,t_0+\alpha R''\bigr],\\
	&h'(t_2)>0,\qquad f'(t_2)=\cos\biggl(\frac{R}{N}\biggr)
\end{align*}
for given $N'',R'',\rho'',N,R>0$, and such that the resulting metric has ${\rm Ric}_k>0$. By replacing~$h(t)$ and $f(t)$ by $\frac{1}{\alpha}h(\alpha(t-t_0))$ and $\frac{1}{\alpha}f(\alpha(t-t_0))$, respectively, for suitable $\alpha>0$, and with an appropriate choice of $N''$ and $R''$, we can achieve $h(t)=\sin(t)$ and $f'(t)=\rho'$ for all $t\in[0,1/2]$ and all $\rho'>0$ small enough, so that \ref{EQ:SCAL_FCTS1}, \ref{EQ:SCAL_FCTS4} and \ref{EQ:SCAL_FCTS5} are automatically satisfied. Further, it follows directly from the construction that \ref{EQ:SCAL_FCTS2} and \ref{EQ:SCAL_FCTS6} hold. With a view towards \ref{EQ:SCAL_FCTS3} we choose $N$, $R$ so that $\cos(R/N)=\Delta$, and hence $f'(t_2)=\Delta$.

Next we show that $t_2$ depends smoothly on $\rho'$. For that, recall that for $t>\varepsilon$, where $\varepsilon>0$ can be chosen arbitrarily small, the functions $h$ and $f$ coincide with the functions $h_a=a\cdot h_0$ and $f_{b,C}=b\cdot f_C$ from \cite{Re}, where $b=\rho'$. These functions are defined on all of $[0,\infty)$ and by \cite[proof of Theorem 4.1]{RW}, for all $b$ sufficiently small, there exists $t_b>0$, so that we have $f'_{b,C}(t_b)>\Delta$ and ${\rm Ric}_k>0$ on $[0,t_b]$. The constant $t_2$ is then defined by $f'_{b,C}(t_2)=\Delta$. Since $f'_{b,C}, f''_{b,C}>0$ on $(0,\infty)$ and $f'_{b,C}(0)=0$, $t_2$ is uniquely determined and depends smoothly on~$b$ and thus on $\rho'$. Since $f_{b,C}=b\cdot f_C$ and $f'_{C}(t)\to\infty$ as $t\to \infty$ by \cite[Lemma 3.7]{Re}, it follows that~$t_2$ increases monotonically to $\infty$ as $\rho'\to 0$. In the following, we consider $f$ and $h$ to be defined on all of $[t_0,\infty)$ by identifying the functions with $f_{b,C}$ and $h_a$ on $[\varepsilon,\infty)$.

Further, differentiating $b\cdot f_C(t_2(b))$, where we write $t_2(b)$ instead of $t_2$ to indicate that $t_2$ depends on $b$, we obtain
\[f_C(t_2(b))+bf_C'(t_2(b))t_2'(b)\geq f_C(t(b))\to\infty \]
as $b\to0$, which establishes \ref{EQ:SCAL_FCTS9} if we set $a=t_2$ (as $h$ is unaffected by $\rho'$). Finally, since changing~$\rho'$ results in scaling the function $f$ and leaving $h$ unchanged, both $f$ and $h$ depend smoothly on~$\rho'$.

We will now modify $f$ and $h$ on $[t_2,\infty)$ so that the remaining conditions \ref{EQ:SCAL_FCTS3}, \ref{EQ:SCAL_FCTS7} and~\ref{EQ:SCAL_FCTS8} are satisfied and so that we have ${\rm Ric}_k>0$ for the extended functions. Recall from \cite[Corol\-lary~2.7]{RW} that for ${\rm Ric}_k>0$ the following inequalities must be satisfied:
\begin{align}
	\label{EQ:RIC_K1}
	&-(k-q)\frac{h''}{h}-q\frac{f''}{f}>0,\\
	\label{EQ:RIC_K2}
	&-\frac{h''}{h}+(k-q-1)\frac{1-{h'}^2}{h^2}-q\frac{f'h'}{fh}>0,\\
	\label{EQ:RIC_K3}
	&(k-q)\frac{1-{h'}^2}{h^2}-q\frac{f'h'}{fh}>0,\\
	\label{EQ:RIC_K4}
	&-\frac{f''}{f}-p\frac{f'h'}{fh}+(k-p-1)\frac{1-{f'}^2}{f^2}>0.
\end{align}

In order to modify the functions on $[t_2,\infty)$, we first create a $C^1$-function as follows. For~${t \ge t_2}$, set $L(t)=f(t_2)+\Delta(t-t_2)$. Clearly, the function
\[\bar{f}:=
\begin{cases}
	f(t) \quad \text{ for } t \le t_2, \\
	L(t) \quad \text{ for } t \ge t_2 \\
\end{cases} \]
is $C^1$ at $t=t_2$ and smooth elsewhere. Assume for the moment that the inequalities \eqref{EQ:RIC_K1}--\eqref{EQ:RIC_K4} hold for all $t\neq t_2$ with $f$ replaced by $\bar{f}$. It then follows from \cite[Lemma 3.1]{RW} that we can smooth $\bar{f}$ in a very small interval $t \in (t_2-\epsilon,t_2+\epsilon)$, whilst keeping all the above inequalities satisfied. This results in the desired function $\tilde{f}.$ It remains, therefore, to justify the assumption that extending $f$ by $L$ in \eqref{EQ:RIC_K1}--\eqref{EQ:RIC_K4} for $t \ge t_2$ results in inequalities which hold for all $t \ge t_2$.

To analyse the curvature properties when replacing $f$ by $\tilde{f}$, we need the following lemma.
\begin{Lemma}\label{f_comparison}
	Let $f_1,f_2\colon[T,\infty)\to (0,\infty)$ be smooth functions such that $f_1(T)=f_2(T)$ and $f_1'(T)=f_2'(T)$. If
	\[\frac{f_1''(t)}{f_1(t)}>\frac{f_2''(t)}{f_2(t)} \]
	for all $t>T$, then
	\[\frac{f_1'(t)}{f_1(t)}>\frac{f_2'(t)}{f_2(t)} \]
	for all $t>T$.
\end{Lemma}
\begin{proof}
	The claim holds if and only if
	\[f_1'(t)f_2(t)-f_2'(t)f_1(t)>0 \]
	for all $t>T$. Further, by the hypothesis, we have equality for $t=T$. Differentiation yields
	\[\bigl(f_1'f_2-f_2'f_1\bigr)'=f_1''f_2-f_2''f_1, \]
	which is positive on $(T,\infty)$.
\end{proof}

\begin{Lemma}
	Assuming the curvature expressions \eqref{EQ:RIC_K1}--\eqref{EQ:RIC_K3} hold at $t>0$ for positive functions~$f$ and $h$ with $f''>0$, $h''<0$, and $f',h'>0$, replacing $f$ by $\tilde{f}$ does not affect the validity of these inequalities for any $t$.
\end{Lemma}

\begin{proof}
	We do not alter the original functions for $t \le t_2-\epsilon$, and by construction, for $t \in (t_2-\epsilon,t_2+\epsilon)$ the inequalities hold with $\tilde{f}$ replacing $f$. At $t \ge t_2+\epsilon$, we have $\tilde{f}=L$, and since $L'' \equiv 0$ it is automatic that \eqref{EQ:RIC_K1} holds for all such $t$.
	
	For \eqref{EQ:RIC_K2} and \eqref{EQ:RIC_K3}, the only terms which are affected by our modification of $f$ are those involving the expression $f'h'/fh.$ But these terms always appear with a negative coefficient, so by Lemma \ref{f_comparison} we see that our modification results in an {\it increase} in value for these terms, and hence an increase for the left-hand side of the inequalities. Thus \eqref{EQ:RIC_K2} and \eqref{EQ:RIC_K3} continue to hold for all $t \ge t_2+\epsilon.$
\end{proof}

\begin{Lemma}\label{inequality_4}
	Suppose that $($in addition to previous assumptions$)$ we have that $h'f$ is monotonically decreasing. Then inequality \eqref{EQ:RIC_K4} holds for all $t>0$ with $f$ replaced by $\tilde{f}$.
\end{Lemma}
\begin{proof}
	Since \eqref{EQ:RIC_K4} holds at $t\le t_2$, and since $f''(t_2)>0$ and $f(t_2)=L(t_2)$, we have
	\[(k-p-1)\frac{1-\Delta^2}{L(t_2)^2}-p\frac{h'(t_2)\Delta}{L(t_2)h(t_2)}>0. \]
	Hence, since $L''\equiv0$, \eqref{EQ:RIC_K4} holds for $f$ replaced by $L$ at $t=t_2$. For $t>t_2$, we observe that \eqref{EQ:RIC_K4} is equivalent to
	\[ 1-\Delta^2>\frac{p\Delta}{k-p-1} \frac{L h'}{h}. \]
	Since $f''>0$, we have $f(t)>L(t)$ for $t>t_2$, and since $h'f$ is monotonically decreasing while $h$ is monotonically increasing, we have
	\[\frac{L(t)h'(t)}{h(t)}<\frac{f(t)h'(t)}{h(t)}<\frac{f(t_2)h'(t_2)}{h(t_2)}=\frac{L(t_2)h'(t_2)}{h(t_2)}<\frac{k-p-1}{p\Delta}\bigl(1-\Delta^2\bigr) \]
	for all $t>t_2$.
\end{proof}

Note that the assumption in Lemma \ref{inequality_4} is satisfied in our case by \cite[proof of Lem\-ma~3.7\,(3)]{Re}.

We have established that having achieved the correct slope $\Delta$, we can straighten the function~$f$ whilst preserving the curvature inequalities \eqref{EQ:RIC_K1}--\eqref{EQ:RIC_K4}.

We now turn our attention to straightening $h$ to a constant function, again preserving conditions \eqref{EQ:RIC_K1}--\eqref{EQ:RIC_K4}, but also ensuring that \ref{EQ:SCAL_FCTS8} is satisfied, i.e.,
\[\sup_{t \in [0,\infty)}|h''(t)|=\sup_{t \in [0,1/2]}|h''(t)|.\]
First we address the issue of bending $h$ in such a way that the curvature conditions continue to hold. Note that we have to loosen inequality \eqref{EQ:RIC_K1} by replacing positivity by non-negativity in order to accommodate such a deformation. We will address the consequences of this at the end.

Recall from before that $h>0$, $h' \in (0,1)$, $h''<0$ and $\lim_{t \to \infty} h'(t)=0.$ It is also easily verified that $h'''>0$. We wish to deform this function at $t \ge t_2$ to a function $\tilde{h}$ so that $\tilde{h}'(t) =0$ for all $t \ge t_3$ for some $t_3>t_2$, $\tilde{h}''(t) \le 0$ and $\tilde{h}'(t) \ge 0$ for all $t$. Given that we are considering $t \ge t_2$, we have $\tilde{f}''/\tilde{f}\equiv 0$, and so the left-hand side of \eqref{EQ:RIC_K1} is simply $-\tilde{h}''/\tilde{h}$, which is non-negative as required.

Note further that the desired deformation will tend to increase the term $\bigl(1-h'^2\bigr)/h^2$, which appears in both \eqref{EQ:RIC_K2} and \eqref{EQ:RIC_K3}. The term $-f'h'/fh$ will also improve provided we can establish the inequality
\[\frac{h'(t)}{h(t)} \ge \frac{\tilde{h}'(t)}{\tilde{h}(t)}\] for all $t$. In particular, this is vital for showing that \eqref{EQ:RIC_K4} continues to hold when $h$ is replaced~by~$\tilde{h}$.

In order to construct $\tilde{h}$, we introduce a $C^1$-function $\phi(t)$ as follows: set $\phi(t)=t$ for all~${t \le t_2}$, we suppose that $\phi(t)$ is constant for all ${t \ge t_3}$, $\phi'(t) \in [0,1]$ for all $t$, and $\phi''(t) \le 0.$ Set~${\tilde{h}(t)=h(\phi(t))}$, so $\tilde{h}'(t)=h'(\phi(t))\phi'(t)$ and $\tilde{h}''(t)=h''(\phi(t))\phi'^2(t)+h'(\phi(t))\phi''(t).$ It is clear that $\tilde{h}$ has all the required properties except possibly for the last one. We claim that a careful choice of $\phi$ will ensure that this last condition also holds.

To this end, set $\lambda_{t_2}=\frac{h''(t_2)}{h'(t_2)}-\frac{h'(t_2)}{h(t_2)}$. It is easily checked that $\lambda_{t_2}<0.$ Define
\[
\psi'(t)=\begin{cases}
	1,\quad & t\leq t_2,\\
	\lambda_{t_2}(t-t_2)+1, & t\in [t_2,t_2-1/\lambda_{t_2}],\\
	0, & t\geq t_2-1/\lambda_{t_2}
\end{cases}
\]
with $\psi(0)=0$. Then $\psi$ is $C^1$, and we can smooth $\psi$ to the desired function $\phi$ in arbitrarily small neighbourhoods of $t=t_2$ and $t=t_2-1/\lambda_{t_2}$, in such a way that $\phi$ is $C^1$-arbitrarily close to $\psi$ and~${\phi''\le 0}$. Moreover, we can arrange for $\phi''$ to lie in the interval $[\lambda_{t_2},0]$, and in a~neighbourhood of $t=t_2$, we will have $\phi'(t)\le \psi'(t)$.

\begin{Lemma}\label{h}
	With $\phi(t)$ as above and $\tilde{h}(t):=h(\phi(t))$, we have
	\[
	\frac{h'(t)}{h(t)} \ge \frac{\tilde{h}'(t)}{\tilde{h}(t)}\]
	for all $t$.
\end{Lemma}

\begin{proof}
	We only need consider $t \ge t_2.$ We have
	\begin{align*}
		\tilde{h}'(t)&{}=h'(\phi(t))\phi'(t).
	\end{align*}
	As $h''<0$, we have that $h'(t_2) \ge h'(\phi(t))$ for all $t \ge t_2$. Similarly, as $h'\ge 0$, we have $h(t_2) \le h(\phi(t))$ for all $t \ge t_2$. Therefore it suffices to show that \[\frac{h'(t)}{h(t)} \ge \frac{h'(t_2)\phi'(t)}{h(t_2)},\] or equivalently
	\begin{gather}\label{ast}
		\phi'(t) \le \frac{h'(t)h(t_2)}{h(t)h'(t_2)}. \tag{$\ast$}
	\end{gather}
	Notice that the right-hand side is equal to 1 at $t=t_2$, its derivative equals $\lambda_{t_2}$ at $t=t_2$, and a~calculation shows that its second derivative is positive. As $\phi'(t_2)\le \psi'(t_2)=1$ and $\phi'' \le 0$, it follows automatically that \eqref{ast} holds.
\end{proof}
It is now clear that all required properties for $\tilde{h}$ defined as in Lemma \ref{h} hold with $t_3=t_2-1/\lambda_{t_2}$ and that inequality \eqref{EQ:RIC_K4} holds for all $t$. We also conclude that both terms on the left-hand side of \eqref{EQ:RIC_K3} improve after this deformation, so \eqref{EQ:RIC_K3} continues to hold as well.

We must be a little careful over \eqref{EQ:RIC_K2}, as this involves a lower coefficient for the $\bigl(1-h'^2\bigr)/h^2$ term than that in \eqref{EQ:RIC_K3}, but also includes a $-h''/h$ term. It is clear from the construction of the functions in \cite{RW} that $h$ can be chosen so that the following inequality (which implies both \eqref{EQ:RIC_K2} and \eqref{EQ:RIC_K3}) holds for all $t$: \[(k-q-1)\frac{1-h'^2}{h^2}-q\frac{f'h'}{fh}>0.\] As the deformation improves the left-hand side, we see that with this particular choice of $h$, \eqref{EQ:RIC_K2} continues to hold when $h$ is replaced by $\tilde{h}.$

\begin{Lemma}
	For $\rho'$ sufficiently small, or, equivalently, for $t_2$ sufficiently large, we have
\[
\sup_{t \in [0,\infty)}\bigl|\tilde{h}''(t)\bigr|=\sup_{t \in [0,1/2]}\bigl|\tilde{h}''(t)\bigr|.
\]
\end{Lemma}

\begin{proof}
	Notice that $|h''|$ is a strictly decreasing function which tends to 0 as $t \to \infty.$ If $t_2$ is sufficiently large, we can assume that $|h''(t_2)|=\bigl|\tilde{h}''(t_2)\bigr|<|h''(1/2)|/2.$ Further, since $h$ is an increasing function and $h'(t)\to 0$ as $t\to\infty$, for $t_2$ sufficiently large, we can assume that
	\[\bigl|h'(t_2)\lambda_{t_2}\bigr|= \bigl|h''(t_2)\bigr|+\frac{h'(t_2)^2}{h(t_2)}<\bigl|h''(1/2)\bigr|/2. \]
	It follows that for $t\in [t_2,t_3]$ we have
	\begin{align*}
		\bigl|\tilde{h}''(t)\bigr|&{}= \bigl|h''(\phi(t))\bigr|\phi'(t)^2+h'(\phi(t))\bigl|\phi''(t)\bigr|
		\leq \bigl|h''(t_2)\bigr|+\bigl|h'(t_2)\lambda_{t_2}\bigr|\\
		&<\frac{\bigl|h''(1/2)\bigr|}{2}+\frac{\bigl|h''(1/2)\bigr|}{2}=\bigl|h''(1/2)\bigr|.\tag*{\qed}
	\end{align*}\renewcommand{\qed}{}
\end{proof}

Note that by construction, $t_3$ depends smoothly on $t_2$ and therefore smoothly on $\rho'$, and setting $a=t_3+\varepsilon$ for some $\varepsilon>0$ so that $\tilde{h}'$ vanishes in a neighbourhood of $a$, yields the required constant $a$. It is also clear that \ref{EQ:SCAL_FCTS9} still holds for this definition of $a$.

Finally, let us return to the issue of non-negativity in inequality \eqref{EQ:RIC_K1} when $f$ and $h$ are replaced by $\tilde{f}$ and $\tilde{h}$. This could lead to a violation of the ${\rm Ric}_k>0$ condition, but it is easy to implement a further small modification which preserves positivity in \eqref{EQ:RIC_K1}. For $t$ slightly less than~$t_3$ \big(where $\tilde{h}'$ hits 0 for the first time\big), introduce a small concave-down bend in the function~$\tilde{f}$. This automatically has the desired effect on \eqref{EQ:RIC_K1}, while boosting the left-hand sides of~\eqref{EQ:RIC_K2}--\eqref{EQ:RIC_K4}. To complete the surgery, we need to perform such a deformation anyway. In \cite{Wr2} this deformation is performed after $h$ is constant, however this is only for the convenience of the exposition: a slight bend $\tilde{f}$ just before $h$ is constant has no implications for curvature-preserving surgery. Moreover, condition \ref{EQ:SCAL_FCTS3}, which now fails to hold after this further deformation of~$\tilde{f}$, merely reflects the approach of \cite{Wr3}. It is easily checked (by consulting~\mbox{\cite[Section~4]{Wr3}}) that this condition has precisely no role to play in the key metric deformation in that paper. Hence making the proposed deformation of $\tilde{f}$ will have no effect on the application of the results of~\cite{Wr3} to our current ${\rm Ric}_k>0$ situation, so we are free to omit that condition from our considerations.

Thus, we constructed functions $\tilde{f}$ and $\tilde{h}$ that satisfy all properties required in \cite[Theorem~3.2]{Wr3} (except \ref{EQ:SCAL_FCTS3} as discussed above), which completes the proof of Theorem~\ref{spaces}.

\subsection*{Acknowledgements}

The first named author acknowledges funding by the SNSF-Project 200020E\textunderscore 193062 and the DFG-Priority programme SPP 2026. Both authors would like to thank Diarmuid Crowley and David Gonz\'alez \'Alvaro for helpful conversations, and the anonymous referees for their careful reading and insightful comments.

\pdfbookmark[1]{References}{ref}
\LastPageEnding

\end{document}